\documentclass[12pt]{amsart}
\usepackage{amsmath,amssymb,amsfonts,amsthm,bbm}
\usepackage[utf8]{inputenc}
\usepackage[T1]{fontenc}
\usepackage{amscd, graphics,url}
\usepackage{color}
\usepackage{hyperref}
\usepackage[margin=2.5cm]{geometry}
\usepackage{dsfont}
\usepackage[mathscr]{euscript}
\usepackage{enumitem}
\numberwithin{equation}{section}
\newtheorem{thm}{Theorem}[section]
\newtheorem{lem}[thm]{Lemma}
\newtheorem{prop}[thm]{Proposition}
\newtheorem{cor}[thm]{Corollary}
\theoremstyle{definition}
\newtheorem{dfn}[thm]{Definition}
\theoremstyle{remark}
\newtheorem{rmk}[thm]{Remark}

\theoremstyle{remark}
\newtheorem{example}[thm]{Example}

\newcommand{\ck}{{\mathcal K}}
\newcommand{\cb}{{\mathcal B}}
\newcommand{\ca}{{\mathcal A}}
\newcommand{\ch}{{\mathcal H}}

\newcommand{\cn}{{\mathcal N}}

\newcommand{\cl}{{\mathcal L}}

\newcommand{\N}{\mathbb{N}}
\newcommand{\Q}{\mathbb{Q}}
\newcommand{\R}{{\mathbb{R}}}

\newcommand{\T}{{\mathbb{T}}}

\newcommand{\Z}{\mathbb{Z}}

 \newcommand{\vw}{{\vec \omega}}
 \newcommand{\vb}{{\vec b}}

\newcommand{\ve}{{\vec e}}

\newcommand{\vn}{{\vec n}}

\newcommand{\vs}{{\vec s}}

\newcommand{\vu}{{\vec u}}
\newcommand{\vv}{{\vec v}}

\usepackage{setspace}
\begin{document}
\title{Directional weak mixing and directional Kronecker for $\Z^2$-actions}
\author{E. Arthur Robinson, Jr.}
\address{Department of Mathematics, George Washington University, Washington, DC 20052, USA.}
\author{ Ay\c{s}e A. \c{S}ahin}
\address{Department of Mathematics, Embry-Riddle Aeronautical University, Daytona Beach, FL 32114-3900, USA.}
\date{October 29, 2026}
\dedicatory{In memory of Joseph Rosenblatt.}
\maketitle

\begin{abstract}
Let $T$ be a measure preserving $\Z^2$-action on a Lebesgue probability space $(X,\mu)$ and let 
$L\subseteq\R^2$ be a $1$-dimensional subspace. For a rational subspace $L$, direction $L$ weak mixing corresponds to 
weak mixing for the measure preserving transformation $T^{\vv_0}$, where $\vv_0$ is a primitive vector in $\Z^2$ generating $L$. 
For an irrational direction $L$, we study two different definitions of direction $L$ weak mixing, and show that they are the same:   
``Dynamical'' direction $L$ weak mixing is defined to be the absence of any direction $L$ Kronecker factors. These factors are defined in 
terms of almost periodicity along a finite width direction $L$ strip in $\Z^2$.   ``Spectral'' direction $L$ weak mixing is defined 
as the absence of any direction $L^\perp$ line segments $I\subseteq\T^2$  with positive spectral measure.
We also give a spectral characterization of the directional Kronecker property defined in \cite{liuxu}, and 
provide examples of $\Z^2$-actions $K$ that are direction $L$ Kronecker. 
\end{abstract}

\bigskip

\section{Introduction}

Let $T$ be an ergodic measure preserving action of $\Z^2$ on a Lebesgue probability space $(X,\mu)$. 
We call a $1$-dimensional  subspace $L\subseteq\R^2$ a \emph{direction} in $\Z^2$.
This paper studies two different 
notions of what it means for $T$  to be \emph{weakly mixing in the direction $L$}. One notion, 
which comes  from  \cite{liu}, and that 
we call here \emph{dynamical direction $L$ weak mixing}, 
is based on the absence of \emph{direction $L$ Kronecker factors}, which 
result from  \emph{direction $L$ almost periodicity}, \cite{liuxu}.
The second notion, which comes from \cite{RRS},
 and that we call here \emph{spectral direction $L$ weak mixing},
 occurs when the spectral measure $\sigma^T$ satisfies $\sigma^T(I)=0$ 
for every line segment $I$ in $\T^2$ perpendicular to $L$.
Our main result, Theorem~\ref{main}, is  that  these two kinds of direction $L$ weak mixing
are, in fact, the same. 
As the final version of our preprint was being prepared, the authors received the preprint \cite{liuetal} which contains 
a proof of the 
same result, obtained  using a different method.
In addition, we give a new spectral characterization of the \emph{direction $L$ Kronecker} property, studied in 
 \cite{liuxu}, and  provide some concrete examples.
 
Directional weak mixing and the directional 
Kronecker property
belong to a growing list of directional properties for $\Z^d$ (and $\R^d$) actions $T$ that have been studied 
in recent years.  The idea of directional dynamics goes back to   
Milnor's \cite{Milnor}, \cite{Milnor1} study of directional entropy, and continues, notably, 
with Boyle and Lind's 
\cite{BL} study of directional expansiveness. Other papers studying directional dynamics 
include \cite{Danrec}, \cite{FeKa}, \cite{JSRec}, 
\cite{liuxu}, \cite{PIsrael}, \cite{Pugh}, \cite{rs}, and most recently \cite{liuetal}.  
Milnor's key observation is that for $e<d$, an  $e$-dimensional direction of a $\Z^d$-action
need not  be restricted
to rank $e$ sub-lattices of $\Z^d$, but rather can often be taken as  
 an arbitrary $e$-dimensional subspace $L\subseteq\R^d$.
Those directions that do  correspond  
to sub-lattices of $\Z^d$ are called \emph{rational} directions. For simplicity, in this paper we stick to the case of 
$\Z^2$-actions and $e=1$ dimensional directions $L$. A few comments on the case of directional weak mixing for 
$\Z^d$-actions in general 
are given in Section~\ref{conclusion}.

\section{Preliminaries}\label{prelim}

Let $(X,\mu)$ be a nonatomic Lebesgue probability space, (with an implicit complete Lebesgue $\sigma$-algebra $\cb$). 
We will write a measurable measure preserving $\Z^2$-action on  
$(X, \mu)$  as $T=\{T^{\vn}\}_{\vn\in\Z^2}$, where $\vn=(n_1,n_2)\in \Z^2$, 
$T^{\vn}=T_1^{n_1}T_2^{n_2}$, and $T_1$ and $T_2$ are a pair of commuting invertible measure 
preserving transformations of  $(X, \mu)$. We will always assume that $T$ is \emph{ergodic}, (i.e., there are 
no non-constant $T$-invariant functions). 
 
Since $(X,\mu)$ is a Lebesgue  space,
the Hilbert space $L^2(X,\mu)$, with inner product   $\langle f,g\rangle:=\int_Xf\overline{g}\,d\mu$ and norm $||f||_2:=\langle f,f\rangle^{1/2}$
is separable.
Let  $U_T=\{U_T^\vn\}_{\vn\in\Z^2}$
be the $T$-induced \emph{Koopman  representation} of $\Z^2$ on  $L^2(X,\mu)$, 
defined $U_T^{\vn}f(x):=f(T^\vn x)$.
We will usually restrict $U_T$ to the closed invariant subspace 
$L^2_0(X,\mu):=\{\text{const.}\}^\perp\subseteq L^2(X,\mu)$ of integral $0$ functions. 
Let $\T^2:=\R^2/\Z^2\approx[0,1)^2$. 
A point $\vw\in\T^2$ is an \emph{``eigenvalue''} for $T$, with \emph{eigenfunction} $f\in L^2_0(X,\mu)$, if 
\begin{equation}\label{eve}
U_T^\vn f=e^{-2\pi i (\vn\cdot\vw)}f,\quad\vn\in\Z^2.
\end{equation}
Since $T$ is ergodic, there is  no eigenvalue $\vec 0$ for any $f\in L^2_0(X,\mu)$.
Moreover, we can assume without loss of generality that $|f(x)|=1$ since $T$ is ergodic.

We say a $\Z^2$-action $T$ is \emph{weakly mixing}, 
(which we call \emph{spectral} weak mixing in this context),
if $U_T$ has no eigenfunctions $f\in L^2_0(X,\mu)$. 
This implies $T$ is ergodic since an invariant function $f\in L^2_0(X,\mu)$ would be an eigenfunction 
for eigenvalue $\vec 0$. 
At the other extreme, we say $\Z^2$-action $T$ is \emph{Kronecker} if the eigenfunctions 
$f\in L^2_0(X,\mu)$ span a dense subspace of $L^2_0(X,\mu)$.

In what follows we will need some spectral theory (see \cite{kt} or \cite{eli}).
The \emph{spectral measure} for $f\in L^2_0(X,\mu)$ 
is the \emph{finite positive Borel measure} $\sigma_f$ on $\T^2$, defined in terms of its Fourier coefficients by
\begin{equation}\label{dsm}
\widehat\sigma_f(\vn)=\int_{\T^2}e^{-2\pi i (\vn\cdot\vu)}\,d\sigma_f(\vu)=\langle U_T^\vn f,f\rangle,
\end{equation}
for $\vn\in \Z^2$. It is easy to see that $\sigma_f=\delta_\vw$, the unit point mass at $\vw\in\T^2$, if and only if $f$ is an eigenfunction 
with eigenvalue $\vw$, \eqref{eve}.

The  unitary  \emph{multiplication representation} $M^\vn$ 
of $\Z^2$ on $L^2(\T^2,\sigma_f)$ is defined  
$M^\vn g(\vu)= e^{-2\pi i (\vn\cdot\vu)}g(\vu)$. 
Let $Z(f)=\overline{\text{span}}\{U_T^\vn f:\vn\in \Z^2\}\subseteq L^2_0(X,\mu)$ be the closed $U_T$-invariant
\emph{cyclic subspace} generated by $f$. 

\begin{prop}[Spectral theorem for cyclic subspaces]\label{simple}
There is a unitary equivalence $\Phi:Z(f)\to L^2(\T^2,\sigma_f)$ so 
that $\Phi\circ U_T^\vn=M^\vn\circ \Phi$, with $\Phi(f)=1$. Moreover, there is an $f_0\in L^2_0(X,\mu)$,
called a function of  \emph{maximal spectral type} 
 so that $\sigma_f\ll \sigma_{f_0}$ for all $f\in L^2_0(X,\mu)$. 
The measure $\sigma_{f_0}$ of \emph{maximal spectral type} is unique up to equivalence.
\end{prop}

For simplicity, we fix a particular $f_0$ of maximal spectral type, and call $\sigma^T:=\sigma_{f_0}$ the \emph{spectral measure} for $T$.

Let $\ch_\text{wm}$ be the set of $f\in L^2_0(X,\mu)$ so that
\begin{equation}\label{wml}
\lim_{N\to\infty}\frac{1}{(N+1)^2}\sum_{\vn\in [0\dots N]^2}|\langle U_T^\vn f, f\rangle|^2
=\lim_{N\to\infty}\frac{1}{(N+1)^2}\sum_{\vn\in [0\dots N]^2}|\widehat\sigma_f(\vn)|^2=0.
\end{equation}
It is clear that if $f\in L^2_0(X,\mu)$ is an eigenfunction, then $f\notin \ch_\text{wm}$.
Conversely, 
Wiener's lemma (see \cite{katznelson}) says that for $f\in L^2_0(X,\mu)$,  the right hand C\'esaro sum in \eqref{wml}
is $0$ if and only if $\sigma_f$ has no atoms. Thus $T$ is weakly mixing if and only if 
$\sigma^T$ has no atoms.

More generally, the
\emph{Lebesgue decomposition} splits $\sigma^T=\sigma^T_\text{d}+\sigma^T_\text{c}$ into 
 \emph{discrete} and 
 \emph{continuous} 
parts, with $\sigma^T_\text{d}\perp\sigma^T_\text{c}$.
Thus by Proposition~\ref{simple} and Wiener's lemma, $\ch_\text{wm}=\ch_\text{c}$, where $\ch_\text{c}:=\{f:\sigma_f\ll\sigma^T_\text{c}\}$, 
and $\ch_\text{d}:=\{f:\sigma_f\ll\sigma^T_\text{d}\}$ satisfies $\ch_\text{d}=\ch^\perp_\text{c}=\ch^\perp_\text{wm}$. 
The resulting $U_T$-invariant orthogonal decomposition 
$L^2_0(X,\mu)=\ch_\text{d}\oplus\ch_\text{c}=\ch_\text{d}\oplus\ch_\text{wm}$, into \emph{discrete} and \emph{continuous} spectrum parts, 
is called \emph{Koopman-von Neumann decomposition}, \cite{kvn}.

A function $f\in L^2_0(X,\mu)$ is called \emph{almost periodic} if $O(f):=\{U_T^\vn f:\vn\in \Z^2\}$ is pre-compact in $L^2_0(X,\mu)$.
Clearly any eigenfunction is almost periodic, and the set $\ch_\text{d}\subseteq \ch_\text{ap}\subseteq L^2_0(X,\mu)$ is a closed 
invariant subspace. But one can show that, in fact, $\ch_\text{d}=\ch_\text{ap}$, providing the 
\emph{Jacobs-de\! Leeuw-Glicksberg decomposition} 
 $L^2_0(X,\mu)=\ch_\text{ap}\oplus\ch_\text{wm}$, \cite{jacobs},\cite{deleeuw}, (see also \cite{moreira}), which is equivalent to 
 the Koopman-von Neumann decomposition. 
 
 In addition to almost periodic functions,
  a set $B\in\cb$ is called \emph{almost periodic} if $O(B):=\{T^{-\vn} B: \vn\in \Z^2\}$ is
 pre-compact in the metric $d(A,B)=\mu(A\triangle B)$ on $\cb$ (equivalent to the $L^2$ norm on characteristic functions). 
 The almost periodic 
 sets define a sub-$\sigma$-algebra $\ck_\text{ap}\subseteq\cb$, which is the smallest sub-$\sigma$-algebra making all 
 $f\in \ch_{\text{ap}}=\ch_{\text{d}}$ measurable. Then $T$ is weakly mixing, 
 (we call it \emph{dynamical} weak mixing in this context), 
 if and only if 
$\ck_\text{ap}$ is trivial. More generally, since $\ck_\text{ap}$ is $T$-invariant  it defines a factor $\Z^2$-action $K$ of $T$, that is Kronecker and  
unique up to isomorphism.
 
The two versions of direction $L$ weak mixing that we study in this paper correspond to a directional version of 
spectral weak mixing (from  \cite{RRS}) and a directional version of 
dynamical weak mixing (from \cite{liu}). Our goal now is to show they are the same.

 \section{Direction $L$ weak mixing}

By a \emph{direction} $L$ in $\Z^2$, we mean a 1-dimensional subspace of $\R^2$. 
We assume (essentially without loss of generality) that $L$ has slope $\beta>0$, and write 
$L=\text{span}\{\vb\}$, where   $\vb:=(1,\beta)$. 
Then $L^\perp=\text{span}\{\vb_\perp\}$ where $\vb_\perp:=(-\beta,1)$. 
Letting $\pi:\R^2\to\T^2=\R^2/\Z^2$, we let $L^\perp_{\vec 0}:=\pi(L^\perp)\subseteq\T^2$, a subgroup, 
and for $\vv\in\T^2$, $L^\perp_\vv:=L^\perp_{\vec 0}+\vv$, a coset. 

Fix  $r>0$ and define 
\[
L(r)=\{(m,n)\in\Z^2: | n-\beta m |<r\},
\]
\[
L^+(r)=\{(m,n)\in\Z^2:0\le m, | n-\beta m |<r\},
\]
and for $\ell\in \N$,
\begin{equation}\label{thick}
L^+(r,\ell)=\{(m,n)\in\Z^2:1\le m\le \ell, | n-\beta m |<r\}.
\end{equation}

\begin{lem}\label{counter}
Suppose $\beta>0$ is irrational and $r\in\N$. Then for $m\in \{1,\dots,\ell\}$, 
we have $(m,n)\in L^+(r,\ell)$ if and only if 
$n\in \{\lfloor \beta m\rfloor-r+1,\ \lfloor \beta m\rfloor-r+2,\ \dots,\ \lfloor \beta m\rfloor+r\}$.
In particular, $\#(L^+(r,\ell))=2r\ell$.
\end{lem}

\subsection{``Dynamical'' direction $L$ weak mixing and direction $L$ Kronecker}
The following definitions come from \cite{liu} and \cite{liuxu}.

For $B\in\cb$, define 
$O_L(B):=\{T^{-\vn} B:\vn\in L(r)\}$, and define 
the $T$ invariant  \emph{direction $L$ Kronecker sub-algebra} 
$\ck_{\text{ap}_L}=\{B\in\cb: O_L(B)\text{ is pre-compact}\}\subseteq\cb$,
in the metric $d(A,B):=\mu(A\triangle B)$. 
It is shown in \cite{liuxu} that $\ck_{\text{ap}_L}$ is $T$-invariant 
and does not depend on $r$.

\begin{dfn}[\cite{liu}]\label{def1}
An ergodic $\Z^2$-action $T$ is said to be \emph{direction $L$ weakly mixing}, in the sense that we call 
\emph{dynamical direction $L$ weakly mixing},  
if $\ck_{\text{ap}_L}=\{X,\emptyset\}$, (i.e., $\ck_{\text{ap}_L}$ is trivial).
\end{dfn}

In the case that $L$ has rational slope $\beta\in\Q$, it is shown \cite{liu} that direction $L$ weak 
mixing is equivalent to (ordinary) weak mixing for the measure preserving transformation
$S=T^\vn$ where $\vn$ generates $L\cap\Z^2\ne\{\vec 0\}$ (also see the proof of Proposition~\ref{JdLG}).

For $f\in L^2_0(X,\mu)$, define $O_L(f):=\{U_T^\vn f:\vn\in L(r)\}$, and define \emph{ direction $L$ almost periodic} functions 
\begin{equation}\label{d:ap}
\ch_{\text{ap}_L}:=\{f\in L^2_0(X,\mu):O_L(f)\text{ is $L^2$ pre-compact}\}.
\end{equation} 
Then $\ch_{\text{ap}_L}=\{f\in L^2_0(X,\mu):f\text{ is }\ck_{\text{ap}_L}\text{ measurable}\}$, and direction $L$ weak mixing is equivalent to 
$\ch_{\text{ap}_L}=\{0\}$.
Also define
\begin{equation}\label{dwml}
\ch_{\text{wm}_L}:=\left\{f\in L^2_0(X,\mu):
\lim_{\ell\to\infty}\frac{1}{\#(L^+(r,\ell))}\sum_{\vn\in L^+(r,\ell)}|\langle U_T^\vn f , f \rangle|^2=0
\right\}.
\end{equation}
Both
$\ch_{\text{ap}_L}$ and   
$\ch_{\text{wm}_L}$ are closed, $U_T$-invariant subspaces that do not depend on $r$. It is shown in \cite{liu} that  
$\ch_{\text{wm}_L}=\ch_{\text{ap}_L}^\perp$, which provides a direction $L$ Jacobs-de\! Leeuw-Glicksberg decomposition
$L^2_0(X,\mu)=\ch_{\text{ap}_L}\oplus\ch_{\text{wm}_L}$.
In particular,  dynamical direction $L$ weak mixing is equivalent to 
\[
\lim_{\ell\to\infty}\frac{1}{\#(L^+(r,\ell))}\sum_{\vn\in L^+(r,\ell)}|\langle U_T^\vn f , f \rangle|^2=0 \text{ for all }f\in L^2_0(X,\mu).
\] 
Additional  characterizations of direction $L$ weak mixing from \cite{liu} are given in 
Corollary~\ref{lrrs} below.

\subsection{``Spectral'' direction $L$ weak mixing}
The definition of this type of weak mixing comes from  \cite{RRS}. 
We begin with a definition and  a lemma from  \cite{RRS}.

\begin{dfn}[\cite{RRS}]
Let $L\subseteq\R^2$ be a direction in $\Z^2$ and let  $\sigma$ 
be a  finite positive Borel measure on $\T^2$. 
We say  a finite or infinite line segment 
$I\subseteq \T^2$, with $I\perp L$ is a direction $L^\perp$
\emph{wall} for $\sigma$ if $\sigma(I)>0$. We say 
$\sigma$ is a direction $L^\perp$ \emph{wall measure} if it is concentrated 
on a finite or countable union $\bigcup_k I_k$ of direction $L^\perp$ walls. 
In particular, this means $\sigma\big((\bigcup_k I_k)^c\big)=0$.
\end{dfn}

In the following generalized Lebesgue decomposition, walls play the role of direction $L^\perp$ ``atoms'' in $\sigma$. 

\begin{lem}[See Proposition~5.6 in \cite{RRS}]\label{rrsleb}
Let $\sigma$ be a  finite positive Borel measure
on $\T^2$ and let 
$L$ be a direction in $\Z^2$.
Then there is a unique decomposition $\sigma=\sigma_{\text{\rm d}_L}+\sigma_{\text{\rm c}_L}$, with  
$\sigma_{\text{\rm d}_L}\perp\sigma_{\text{\rm c}_L}$, such that
$\sigma_{\text{\rm d}_L}$ is a direction $L^\perp$ wall measure, and 
$\sigma_{\text{\rm c}_L}(I)=0$ has no direction $L^\perp$ walls {\rm (}i.e., $\sigma_{\text{\rm c}_L}(I)=0$ for all $I\perp L${\rm )}.
\end{lem}

We will apply Lemma~\ref{rrsleb} to the spectral measure $\sigma^T$
to obtain $\sigma^T=\sigma^T_{\text{\rm d}_L}+\sigma^T_{\text{\rm c}_L}$. 

\begin{dfn}[\cite{RRS}]\label{dlwmrrs}
An ergodic $\Z^2$-action $T$ on a Lebesgue probability space $(X,\mu)$ is said to be \emph{direction $L$ weakly mixing},
in the sense that we call  \emph{spectral}  direction $L$ weak mixing, 
if its spectral measure satisfies $\sigma^T=\sigma^T_{\text{c}_L}$. 
\end{dfn}

More generally, we have the following.

\begin{lem}\label{kvn1}
Define $\ch_{\text{\rm d}_L}:=\{f\in L^2_0(X,\mu): \sigma_f\ll\sigma^T_{\text{\rm d}_L}\}$ and 
$\ch_{\text{\rm c}_L}:=\{f\in L^2_0(X,\mu): \sigma_f\ll\sigma^T_{\text{\rm c}_L}\}$,
both closed and $U_T$-invariant. Then 
$\ch_{\text{\rm d}_L}\perp \ch_{\text{\rm c}_L}$ and $L^2_0(X,\mu)=\ch_{\text{\rm d}_L}\oplus\ch_{\text{\rm c}_L}$.
Spectral direction $L$ weak mixing is equivalent to $\ch_{\text{\rm c}_L}=L^2_0(X,\mu)$. 
\end{lem}

This follows from Lemma~\ref{rrsleb} and Definition~\ref{dlwmrrs}.

\begin{rmk}\label{fixdef}
 Definition~\ref{dlwmrrs} is based on Proposition 4.5 in \cite{RRS}. 
The original definition of direction $L$ weak mixing  in \cite{RRS} , which is equivalent to this one, 
is stated in terms of the unit-suspension $\R^2$ 
action $\widetilde T$ on $(X\times[0,1)^2,\mu\times d\vs\,)$. That definition is that the restriction of $\widetilde T$ to the direction $L$ 
(essentially an $\R$-action)
has no eigenfunctions orthogonal to the  functions that depend only on $[0,1)^2$. 
A similar remark applies to Definition~\ref{dlwmrrs2} below. 
\end{rmk}

In the case that $L$ has rational slope $\beta$, it is shown in \cite{RRS} that, like dynamical direction $L$ weak mixing,
spectral direction $L$ weak mixing for $T$
is equivalent to (ordinary) weak mixing for the measure preserving transformation
$S=T^\vn$ where $\vn$ generates $L\cap\Z^2\ne\{\vec 0\}$ (also see the proof of Proposition~\ref{JdLG}). 

\subsection{The main theorem}\label{s:main}

\begin{thm}\label{main}
Dynamical and spectral direction $L$ weak mixing are the same.
\end{thm}

This is a corollary of the following more general proposition.
 
\begin{prop}\label{JdLG} 
Let $T$ be an ergodic $\Z^2$-action  on a Lebesgue probability space $(X,\mu)$.
Then $\ch_{\text{\rm d}_L}=\ch_{\text{\rm ap}_L}$ and $\ch_{\text{\rm c}_L}=\ch_{\text{\rm wm}_L}$.
In particular, this provides direction $L$  
\emph{Koopman-von Neumann} and \emph{Jacobs-de Leeuw-Glicksberg}  decompositions 
$L^2_0(X,\mu)=\ch_{\text{d}_L}\oplus\ch_{\text{c}_L}=\ch_{\text{ap}_L}\oplus\ch_{\text{wm}_L}$.
\end{prop}

\section{A directional Wiener lemma}\label{s:wl}

Our main result in this section is the following direction $L$ version of Wiener's lemma on atoms.

\begin{thm}\label{wiener}
Let $\sigma$ be a 
finite positive Borel measure on $\T^2$, where 
$\sigma=\sigma_{\text{\rm d}_L}+\sigma_{\text{\rm c}_L}$
as in Lemma~\ref{rrsleb}.  Let
\begin{equation}\label{ft}
\widehat\sigma(m,n):=\int_{\T^2}e^{-2\pi i (mu+nv)}\,d\sigma(u,v),
\end{equation}
for $(m,n)\in\Z^2$,
be the Fourier-Stieltjes coefficients of $\sigma$. Suppose
$L$ is a direction with irrational slope $\beta>0$. For $r,\ell\in \N$, define  
\begin{equation}\label{al}
W_\ell:=
\frac{1}{\#(L^+(r,\ell))}
\sum_{(m,n)\in L^+(r,\ell)}
|\widehat\sigma(m,n)|^2.
\end{equation}
Then
\begin{enumerate}
\item
$\sigma_{\text{d}_L}=0$ implies
$\lim_{\ell\to\infty} W_\ell = 0$, and 
\item  $\sigma_{\text{d}_L}\ne 0$ implies 
$\lim_{\ell\to\infty} W_\ell > 0$, (i.e., the limit exists, is real and positive).
\end{enumerate}
\end{thm}

The proof occupies the remainder of this section, starting with some preliminary results. 
  We begin with the following, which is a special case of Theorem~\ref{wiener}.

\begin{prop} \label{single}
Let $L$ be a direction with irrational slope $\beta>0$.
Suppose $\sigma_\theta$ 
is a positive finite  Borel measure on $\T^2$ concentrated on 
$L^\perp_{\theta\ve_1}$. 
Define
a measurable bijection 
\begin{equation}
\psi_\theta:\R\to L^\perp_{\theta\ve_1},\quad\quad 
\psi_\theta(t):=\{t\vb_\perp+\theta\ve_1\},
\end{equation} 
where $\vu\mapsto \{\vu\}:=\vu-\lfloor \vu \rfloor\in [0,1)^2$.
Let  $\nu_\theta:=\sigma_\theta\circ\psi_\theta$ be the 
pull-back of $\sigma_\theta$ by $\psi_\theta$, a finite Borel measure on $\R$.

Then 
\begin{equation}\label{dirdir}
\lim_{\ell\to\infty}\frac{1}{\# L^+(r,\ell)}\sum_{(m,n)\in L^+(r,\ell)}|\widehat \sigma_\theta(m,n)|^2
=\frac{1}{2r}\int_{-r}^r|\widehat \nu_\theta(\omega)|^2\,d\omega>0,
\end{equation}
where $\widehat\nu_\theta(\omega):=\int_\R e^{-2\pi i \omega t}\,d\nu_\theta(t)$ is the Fourier transform of $\nu_\theta$.
\end{prop}

\begin{proof}
By the change of variables formula, the Fourier coefficients \eqref{ft} of $\sigma_\theta$ satisfy 
\begin{align*}
\widehat\sigma_\theta(m,n)&:=\int_{L^\perp_{\theta\ve_1}}e^{-2\pi i ((m,n)\cdot(u,v))}d\sigma_\theta(u,v)
=\int_\R e^{-2\pi i ((m,n)\cdot\psi_\theta(t))}d\nu_\theta(t)\\
&=\int_\R e^{-2\pi i ((n-\beta m)t+\theta m)}d\nu_\theta(t)=e^{-2\pi i m\theta}\int_\R e^{-2\pi i (n-m\beta)t}d\nu_\theta\\
&= e^{-2\pi i m\theta}\,\widehat\nu_\theta(n-m\beta),
\end{align*}
 and thus
 \begin{equation}\label{dr}
 |\widehat\sigma_\theta(m,n)|^2=|\widehat\nu_\theta(n-m\beta)|^2.
 \end{equation}

Let 
\begin{equation}\label{irr}
R:\T\to\T,\quad R^m(x):=\{x+\beta m\}
\end{equation}
be the \emph{irrational rotation} by $\beta$ on the circle $\T=\R/\Z\approx[0,1)$.
We view $R$ as a \emph{uniquely ergodic homeomorphism} .
Then $R^m(0)=\{\beta m\}$ 
 so that 
$\beta m=\lfloor\beta m\rfloor +R^m(0)$, and  
\[
n-\beta m=n-\lfloor\beta m\rfloor-R^m(0)=k-R^m(0),\text{ where }k=n-\lfloor\beta m\rfloor.
\]
By \eqref{thick} and \eqref{dr}, 
\begin{align}
\sum_{(m,n)\in L^+(r,\ell)}|\widehat\sigma_\theta(m,n)|^2&=
\sum_{m=1}^\ell\sum_{\overset{n\in\Z}{|n-\beta m|<r}}|\widehat\nu_\theta(n-m\beta)|^2\nonumber \\
&=\sum_{m=1}^\ell\sum_{\overset{k\in\Z}{|k-R^m(0)|<r}}|\widehat\nu_\theta(k-R^m(0))|^2
=\sum_{m=1}^\ell h(R^m(0))\label{stp1},
\end{align}
where
\begin{equation}\label{defh}
h(x):=\sum_{\overset{k\in\Z}{|k-x|<r}}|\widehat \nu_\theta(k-x)|^2
=\sum_{\overset{k\in\Z}{|x-k|<r}}|\widehat \nu_\theta(x-k)|^2
=\sum_{k\in\Z}\mathbbm{1}_{(-r,r)}(x-k)|\widehat\nu_\theta(x-k)|^2,
\end{equation}
since $|x-k|=|k-x|$, and because $\nu_\theta$ is a finite real Borel measure,
$\widehat\nu_\theta(-\omega)=\overline{\widehat\nu_\theta(\omega)}$.

Since $h(x)$ is piecewise continuous, it is Riemann integrable, 
so the unique ergodicity of the irrational rotation $R$ implies  
\begin{equation}\label{stp3}
\lim_{\ell\to\infty}\frac{1}{\ell}\sum_{m=1}^\ell h(R^m(0))=\int_0^1 h(x)\,dx.
\end{equation}
Thus
\begin{align*}
\lim_{\ell\to\infty}\frac{1}{\# L^+(r,\ell)}\sum_{L^+(r,\ell)}|\widehat\sigma_\theta(m,n)|^2
&= \frac{1}{2r}\left(\lim_{\ell\to\infty}\frac{1}{\ell}\sum_{m=1}^\ell h(R^m(0))\right)\\
&=\frac{1}{2r}\int_0^1h(x)\,dx,
\end{align*}
by Lemma~\ref{counter}, \eqref{stp1} and \eqref{stp3}.

Using Tonelli's theorem, and then changing variables $\omega=x-k$ for each $k\in\Z$,
\begin{align*}
\int_0^1 h(x)\,dx&
=\sum_{k\in\Z} \int_0^1 \mathbbm{1}_{(-r,r)}(x-k)|\widehat\nu_\theta(x-k)|^2\,dx
=\sum_{k\in\Z}\int_{-k}^{-k+1} \mathbbm{1}_{(-r,r)}(\omega)|\widehat\nu_\theta(\omega)|^2\,d\omega\\
&=\int_\R \mathbbm{1}_{(-r,r)}(\omega)|\widehat\nu_\theta(\omega)|^2\,d\omega
=\int_{-r}^r |\widehat\nu_\theta(\omega)|^2\,d\omega.
\end{align*}
This integral is positive since $|\widehat\nu_\theta(\omega)|^2$ is  continuous, non-negative, and
$|\widehat\nu_\theta(0)|^2=
\nu_\theta(\R)^2>0$. 
\end{proof}

The next result is a direction $L$ version of the one-sided C\'esaro sum of a Dirichlet kernel.

\begin{prop}\label{limlem}
Let $L$ be a direction with irrational slope $\beta>0$, and let
$r\in\mathbb N$. For $(u,v)\in\mathbb T^2=\mathbb R^2/\mathbb Z^2$, let
\begin{equation}\label{wl}
A_\ell(u,v):=
\frac{1}{\#(L^+(r,\ell))}\sum_{(m,n)\in L^+(r,\ell)} e^{-2\pi i (m u + n v)}
\end{equation}
Then the limit 
\begin{equation}\label{Ainf}
A_\infty(u,v):=\lim_{\ell\to\infty}A_\ell(u,v)
\end{equation} 
exists {\rm (}pointwise{\rm )} for all $(u,v)\in\T^2$, and 
$A_\infty(u,v)=0$
for $(u,v)\notin L^\perp_{\vec 0}$.
\end{prop}

For the proof we use a version of Theorem 1.1 from \cite{rww}, a uniformly convergent Wiener-Wintner theorem.

\begin{thm}[\cite{rww}]\label{wwt}
 Let $R$ be a uniquely ergodic
homeomorphism of a compact metric space $X$. Assume that for every ``eigenvalue'' $w\in \T\approx [0,1)$ 
the corresponding eigenfunction $f_w$ {\rm (}i.e. $f_w(Rx)=e^{-2\pi i w}f_w(x)${\rm )} can be chosen continuous. 
If $g(t)$ a Riemann integrable function\footnote{
The result in \cite{rww} applies to continuous $g$, but extends to Riemann integrable functions $g$ by a standard argument.
}
on $X$, then the limit
\begin{equation}\label{wweq}
\lim_{\ell\to\infty}\frac{1}{\ell}\sum_{m=1}^{\ell} g(R^m (x))e^{-2\pi i w m}
\end{equation}
converges uniformly for $x\in X$. If $w$ is not an eigenvalue then the convergence is to $0$.
In particular, convergence of \eqref{wweq} holds when  $R$ is the irrational rotation by $\beta$ on the circle $\T$. In this case, 
the value of \eqref{wweq} 
is $0$ whenever
$w\ne R^m(0)$ for any $m\in\Z$.
\end{thm}

For the irrational rotation $R$, the eigenfunctions are $f_n(x)=e^{2\pi i n x}$, $n\in\Z$. These satisfy
$f_n(Rx)=e^{2\pi i \beta n}f_n(x)=e^{-2\pi i w}f_n(x)$ for eigenvalue $w=\{-\beta n\}=R^{-n}(0)$, and in particular
the set of eigenvalues is symmetric under negation.

\begin{proof}[Proof of Proposition \ref{limlem}]
By  Lemma~\ref{counter},
\begin{equation}\label{work}
\sum_{\substack{n\in\mathbb Z\\ |n-\beta m|<r}} e^{-2\pi i n v}
=e^{-2\pi i \lfloor \beta m\rfloor v}\sum_{j=-r+1}^{r} e^{-2\pi i j v}
=2r\,D_r(v)\,e^{-2\pi i \lfloor \beta m\rfloor v},
\end{equation}
where
\begin{equation*}
D_r(v):=\frac{1}{2r}\sum_{j=-r+1}^{r}e^{-2\pi i j v}
=\frac{e^{-\pi i v}}{2r}\,\frac{\sin(2\pi r v)}{\sin(\pi v)}
\end{equation*}
is a 2-sided \emph{Dirichlet kernel}. 
Let 
\begin{equation}\label{theta}
\theta(u,v):=\{u+\beta v\}\in[0,1)
\end{equation}
Starting with \eqref{wl}, we apply \eqref{work},  $\lfloor \beta m\rfloor=\beta m-\{\beta m\}$, and  
\eqref{theta}, to obtain 
\begin{align}
A_\ell(u,v)&=D_r(v)\cdot \frac{1}{\ell}\sum_{m=1}^{\ell} e^{-2\pi i (m u + \lfloor \beta m\rfloor v)}
=D_r(v)\cdot \frac{1}{\ell}\sum_{m=1}^{\ell}e^{-2\pi i (m\theta(u,v) - v\,\{\beta m\})}\nonumber\\
&=D_r(v)\cdot \frac{1}{\ell}\sum_{m=1}^{\ell} e^{-2\pi i m\theta(u,v)}\, g_v(R^m(0))\label{pwl}
\end{align}
Here, 
$g_v(t):=e^{2\pi i v t}$, piecewise continuous, so Riemann integrable, and 
$R$ is the irrational rotation by $\beta$.
From the uniform convergence of Theorem~\ref{wwt} we see that \eqref{wweq} converges at $x=0$.
Thus the limit \eqref{Ainf}  converges for all $(u,v)\in\T^2$.

It remains to determine where the limit \eqref{Ainf} is $0$. 
Let $P$ be the parallelogram spanned by the vectors $\ve_1=(1,0)$ and $\vb_\perp=(-\beta,1)$,
with the top and right sides of $P$ open. Then $P$ is a
\emph{fundamental domain} for $\T^2=\R^2/\Z^2$ (i.e., $\R^2=\cup_{\vn\in\Z^2}(P+\vn)$, where the union is disjoint).

Let $I_\theta$ be the interval of slope $-1/\beta$ that goes from the top to the bottom of $P$, through $(\theta,0)$.
Then for $(u,v)\in P\approx\T^2$ we have  $(u,v)\in I_{\theta(u,v)}$ where $\theta(u,v)=\{u+\beta v\}\in [0,1)$.
Viewed as a subset of $P$, $L^\perp_{\theta\ve_1}=\bigcup_{n\in\Z} I_{\{\theta+\beta n\}}$, and in 
particular, $L^\perp_{\vec 0}=\bigcup_{n\in\Z} I_{\{\beta n\}}$. 
By Theorem~\ref{wwt} the limit is $0$ whenever $w=\theta(u,v)$ 
is \emph{not} an eigenvalue of $R$,
which means $\theta(u,v)\ne\{m\beta\}=R^m(0)$ for any $m\in \Z$. This
is equivalent to $(u,v)\notin L^\perp_{\vec 0}$.
\end{proof}
 
 For a finite Borel measure 
 $\eta$ on $\T^2$,  let 
$\tilde{\eta}(B):=\eta(-B)$. For two such measures, the \emph{convolution} 
$\eta_1*\tilde\eta_2=(\eta_1\times\eta_2)\circ D^{-1}$,
the push-forward 
by $D:\T^2\times\T^2\to\T^2$ defined
$D(\vu_1,\vu_2):=\vu_1-\vu_2$. Thus 
\begin{equation}\label{fc}
\int_{\T^2} f(\vu)\,d(\eta_1*\tilde\eta_2)(\vu)
=\int_{\T^2\times \T^2}f(\vu_1-\vu_2)\, d(\eta_1\times \eta_2).
\end{equation}
for any $f\in C(\T^2)$. 
The \emph{autocorrelation} of $\eta$ is the finite
positive Borel measure on $\T^2$
defined by $\eta*\tilde\eta$, and a calculation shows \begin{equation}\label{tau}
\widehat{(\eta*\tilde\eta)}(m,n)=
|\widehat\eta(m,n)|^2\ge 0.
\end{equation} 

By the arguments at the and of the proof of Lemma~\ref{limlem}, we have 
$L^\perp_{\theta\ve_1}=L^\perp_{\theta'\ve_1}$ if and only if $\theta'=R^m(\theta)=\{\theta+\beta m\}$ for some $m\in \Z$, where $R$ is the irrational rotation by $\beta$. We call a countable $\Theta\subseteq\T\approx [0,1)$ a \emph{transversal} if $0\in \Theta$, and 
whenever $\theta,\theta'\in \Theta$ are distinct, 
then $\theta'\ne R^m(\theta)$ for any $m\in \Z$, (i.e., the cosets $L^\perp_{\theta\ve_1}$, $\theta\in \Theta$, are pairwise disjoint). 

We say a finite positive Borel measure $\sigma$ on $\T^2$ is concentrated on a set $E\subseteq\T^2$ if $\sigma(E^c)=0$ 
The next lemma is a refinement of Lemma~\ref{rrsleb}.

\begin{lem}\label{moreleb}
 Let $L\subseteq\R^2$ be a direction
with irrational slope $\beta>0$.
For $\sigma$ a finite positive Borel measure on $\T^2$, suppose 
$\sigma_{\text{d}_L}\ne 0$, where $\sigma=\sigma_{\text{c}_L}+\sigma_{\text{d}_L}$ is the decomposition in 
{\rm Lemma~\ref{rrsleb}}. Then 
$\sigma_{\text{d}_L}$ is concentrated on $\bigcup_{\theta\in \Theta}L^\perp_{\theta\ve_1}$, with  $\Theta$ countable transversal so that 
 $\sigma_{\text{d}_L}(L^\perp_{\theta\ve_1})> 0$ for each $\theta\in \Theta$.
We write 
\begin{equation}\label{deco}
\sigma=\sigma_{\text{c}_L}+\sigma_{\text{d}_L}=\sigma_{\text{c}_L}+\sum_{\theta\in \Theta}\sigma_\theta,
\end{equation}
where each 
$\sigma_\theta$ is the part of $\sigma_{\text{d}_L}$ concentrated on $L^\perp_{\theta\ve_1}$. 
\end{lem}
 
We have
\begin{align}
\tau=\sigma*\tilde\sigma
&=(\sigma_{\text{d}_L}*\tilde{\sigma}_{\text{c}_L}) + 
(\sigma_{\text{c}_L}*\tilde{\sigma}_{\text{d}_L}) +
(\sigma_{\text{c}_L}*\tilde{\sigma}_{\text{c}_L})+
\sum_{\theta_1,\theta_2\in \Theta}
(\sigma_{\theta_1}*\tilde{\sigma}_{\theta_2})\nonumber\\
&=\tau_{{\rm c}_L}+\tau_{{\rm d}_L}=\tau_{{\rm c}_L}+\tau_{L^\perp_{\vec 0}}+\tau_{L^\perp_{\ne \vec 0}},\label{sdlz}
\end{align}
where $\tau_{{\rm c}_L}=(\sigma_{\text{d}_L}*\tilde{\sigma}_{\text{c}_L})+(\sigma_{\text{c}_L}*\tilde{\sigma}_{\text{d}_L})+(\sigma_{\text{c}_L}*\tilde{\sigma}_{\text{c}_L})$, and where
\[
\tau_{L^\perp_{\vec 0}}=\sum_{\theta\in \Theta}(\sigma_{\theta}*\tilde{\sigma}_{\theta}),
\quad
\tau_{L^\perp_{\ne \vec 0}}=\sum_{\overset
{\theta,\theta'\in \Theta}{\theta\ne\theta'}}
(\sigma_{\theta}*\tilde{\sigma}_{\theta'}),
\text{\ \  and\ \  }
\tau_{{\rm d}_L}=\tau_{L^\perp_{\vec 0}}+\tau_{L^\perp_{\ne {\vec 0}}}.
\]

\begin{lem}\label{cutup}
$\tau_{{\rm c}_L}(L^\perp_{\vec 0})=\tau_{L^\perp_{\ne {\vec 0}}}(L^\perp_{\vec 0})=0$, but 
$\tau_{L^\perp_{\vec 0}}$ is concentrated on $L^\perp_{\vec 0}$, and in particular, $\tau_{L^\perp_{\vec 0}}({L^\perp_{\vec 0}})>0$.
\end{lem}

\begin{proof}
For $\eta_1$ and $\eta_2$ finite positive Borel measures
on $\T^2$,
the definition of 
$\eta_1*\widetilde\eta_2$ applied to the subgroup 
$L^\perp_{\vec 0}$ gives
\begin{equation}\label{sbgp}
(\eta_1*\widetilde\eta_2)(L^\perp_{\vec 0})=\int_{\T^2}\eta_1(L^\perp_{\vec 0}+\vu)d\eta_2(\vu),
\end{equation}
$L^\perp_{\vec 0}+\vu=L^\perp_{\theta\ve_1}$ for some $\theta$.
Since $\sigma_{\text{c}_L}$ vanishes on every wall,
$(\sigma_{\text{c}_L}*\tilde{\sigma}_{\text{d}_L})(L^\perp_{\theta \ve_1})=
(\sigma_{\text{c}_L}*\tilde{\sigma}_{\text{c}_L})(L^\perp_{\theta \ve_1})=0$
for all $\theta$.
But also $(\sigma_{\text{d}_L}*\tilde{\sigma}_{\text{c}_L})(L^\perp_{\theta \ve_1})=0$ since $\sigma_{\text{d}_L}$ vanishes off 
$\bigcup_{\theta\in \Theta}L^\perp_{\theta\ve_1}$, which is a null set for  $\sigma_{\text{c}_L}$.
By \eqref{sbgp}, $\sigma_{\theta}*\tilde{\sigma}_{\theta'}$ is concentrated on 
$L^\perp_{(\theta-\theta')\ve_1}$. 
If $\theta,\theta'\in \Theta$ with  $\theta\ne\theta'$, then since $\Theta$ is a transversal, 
$L^\perp_{(\theta-\theta')\ve_1}\ne L^\perp_{\vec 0}$. 
Thus $L^\perp_{\vec 0}\cap\left( \cup_{\overset
{\theta,\theta'\in \Theta}{\theta\ne\theta'}}L^\perp_{(\theta-\theta')\ve_1}\right)=\emptyset$,
which implies $\tau_{L^\perp_{\ne {\vec 0}}}(L^\perp_{\vec 0})=0$.
Finally,  $\tau_{L^\perp_{\vec 0}}=\sum_{\theta\in \Theta}(\sigma_{\theta}*\tilde{\sigma}_{\theta})$ is concentrated on $L^\perp_{\vec 0}$.
Since $\tau=\tau_{{\rm c}_L}+\tau_{L^\perp_{\vec 0}}+\tau_{L^\perp_{\ne \vec 0}}$ is a positive measure,
 and since $\tau_{{\rm c}_L}(L^\perp_{\vec 0})=\tau_{L^\perp_{\ne {\vec 0}}}(L^\perp_{\vec 0})=0$,
  it follows that  $\tau_{L^\perp_{\vec 0}}(L^\perp_{\vec 0})>0$.
\end{proof}

\begin{proof}[Proof of Theorem~\ref{wiener}]
Let $\tau=\sigma*\tilde\sigma$. 
By  \eqref{al} and \eqref{wl},  $W_\ell=\int_{\T^2} A_\ell(u,v)d\tau$.  
It follows from the bounded convergence theorem, Proposition~\ref{limlem}  and Lemma~\ref{cutup} that
 \begin{align*}
 \lim_{\ell\to\infty}W_\ell&= \lim_{\ell\to\infty}\int_{\T^2}A_\ell(u,v)\,d\tau(u,v)\\
 &=\int_{\T^2}A_\infty(u,v)\,d\tau(u,v)\\
 &=\int_{L^\perp_{\vec 0}}A_\infty(u,v)\,d\tau_{L^\perp_{\vec 0}}(u,v),
 \end{align*}

For ({\it 1}\,), if ${\sigma}_{\text{d}_L}=0$ then \eqref{sdlz} shows 
$\tau=\tau_{\text{c}_L}={\sigma}_{\text{c}_L}*{{\tilde\sigma}_{\text{c}_L}}$.  
So $\tau_{\text{d}_L}=0$, which implies $\tau_{L^\perp_{\vec 0}}=0$. Thus, $\lim_{\ell\to\infty}W_\ell=0$.

For ({\it 2}\,), suppose 
$\sigma_{\text{d}_L}\ne 0$, and using  \eqref{deco},
let $\sigma_{\text{d}_L}=\sum_{\theta\in \Theta}\sigma_\theta$,
where each $\sigma_\theta\ne 0$. Let $W_\ell^\theta$ be $W_\ell$ in
\eqref{al} for $\sigma=\sigma_\theta$. 
For each $\theta\in \Theta$, let $\nu_\theta$ be the pull-back to $\R$ of $\sigma_\theta$ from $\T^2$ by $\psi_\theta$,  (as in Proposition~\ref{single}). 
Then
\begin{align*}
\lim_{\ell\to\infty}W_\ell
&=\int_{\T^2}A_\infty(u,v)\,d\tau_{L^\perp_{\vec 0}}
=\int_{\T^2}A_\infty(u,v)\,d\!\left(\sum_{\theta\in \Theta}(\sigma_{\theta}*\tilde{\sigma}_{\theta})\right)
=\sum_{\theta\in \Theta}\left(\int_{\T^2}A_\infty(u,v)\,d(\sigma_{\theta}*\tilde{\sigma}_{\theta})\right)\\
&=\sum_{\theta\in \Theta}\left(\lim_{\ell\to\infty} W_\ell^\theta\right)
=\sum_{\theta\in \Theta}\left(\frac{1}{2r}\int_{-r}^r|\widehat \nu_\theta(\omega)|^2\,d\omega\right)>0.
\end{align*}
\end{proof}

\begin{proof}[Proof of Proposition~\ref{JdLG}] 
First suppose $\beta$ is irrational. We claim that $\ch_{\text{wm}_L}=\ch_{\text{d}_L}^\perp$.  
Since it is shown in \cite{liu} that $\ch_{\text{ap}_L}^\perp=\ch_{\text{wm}_L}$, Lemma~\ref{kvn1} 
then shows $\ch_{\text{ap}_L}=\ch_{\text{d}_L}$ and 
$\ch_{\text{c}_L}=\ch_{\text{wm}_L}$, 
proving Proposition~\ref{JdLG}.

To prove the claim, we first observe that for $f\in L^2_0(X,\mu)$,  
\begin{equation}\label{zero1}
\lim_{\ell\to\infty}\frac{1}{\#L^+(r,\ell)}
\sum_{\vn\in L^+(r,\ell)}|\langle U_T^\vn f,f\rangle|^2
=
\lim_{\ell\to\infty}\frac{1}{\#L^+(r,\ell)}
\sum_{\vn\in L^+(r,\ell)}
|\widehat\sigma_f(\vn)|^2
=
\lim_{\ell\to\infty}W_\ell,
\end{equation}
where $W_\ell$ in \eqref{al} corresponds to $\sigma=\sigma_f$.
Here we have fixed $r\in\N$ without loss of generality.

If $f\in \ch_{\text{d}_L}^\perp$ then by Lemma~\ref{kvn1}, $f\in \ch_{\text{c}_L}$ 
so that $\sigma_f\ll\sigma^T_{\text{c}_L}$, 
and thus by Lemma~\ref{rrsleb}, 
 $\sigma_f(L^\perp_{\theta \ve_1})=0$ for all $\theta\in\T$.
By part (1) of Theorem~\ref{wiener}, $\lim_{\ell\to\infty}W_\ell=0$, 
which implies $f\in\ch_{\text{wm}_L}$.

Conversely, if $f\in\ch_{\text{wm}_L}$ then \eqref{dwml} and \eqref{zero1} 
imply $\lim_{\ell\to\infty}W_\ell=0$ for $\sigma=\sigma_f$.
By part (2) of  
Theorem~\ref{wiener}, this implies that $\sigma_f(L^\perp_{\theta \ve_1})=0$ for all $\theta\in \T$.
Thus $\sigma_f\ll\sigma^T_{\text{\rm c}_L}$ and $f\in \ch_{\text{c}_L}=\ch_{\text{d}_L}^\perp$.

Next we suppose $L$ is a rational subspace. Suppose $L=\text{span}(\vv)$ where $\vv=(v_1,v_2)$ is \emph{primitive}, 
meaning $\text{gcd}(v_1,v_2)=1$. Let $S=T^\vv$, which is a measure preserving transformation of $(X,\mu)$. We have already 
noted that direction $L$ weak mixing is equivalent to the weak mixing of $S$. Our goal is to generalize this result to Proposition~\ref{JdLG}.

Let $(U_S f)(x)=f(Sx)$ be the Koopman operator for $S$ on $L^2_0(X,\mu)$. It is then classical 
that we have orthogonal $U_S$-invariant Koopman-von Neumann 
and Jacobs-de Leeuw-Glicksberg decompositions,
$L^2_0(X,\mu)=\ch^S_{\text{d}}\oplus\ch^S_{\text{c}}=\ch^S_{\text{ap}}\oplus\ch^S_{\text{wm}}$,
where $\ch^S_{\text{d}}=\ch^S_{\text{ap}}$ and $\ch^S_{\text{c}}=\ch^S_{\text{wm}}$.
It is easy to verify that  $\ch^S_{\text{d}}=\ch_{\text{d}_L}$, $\ch^S_{\text{c}}=\ch_{\text{c}_L}$,
$\ch^S_{\text{ap}}=\ch_{\text{ap}_L}$, and $\ch^S_{\text{wm}}=\ch_{\text{wm}_L}$, and the result 
for  $L$ a rational subspace follows.
\end{proof}

\section{Directional weak mixing for $\Z^2$-actions}\label{together}

Having now proved that for ergodic $\Z^2$-actions, the two types of direction $L$ weak mixing  
are the same (Theorem~\ref{main}), 
we simply refer to direction $L$ weak mixing from now on. 
In the following corollaries, we summarize some of the main properties of direction $L$ 
weak mixing which follow from this paper, \cite{liu} and \cite{RRS}. First we mention the following:

\begin{dfn}[\cite{RRS}]\label{dlwmrrs2}
We say an ergodic  $\Z^2$-action is (spectrally)  \emph{direction $L$ ergodic} 
if its spectral measure satisfies $\sigma^T(L^\perp_{\vec 0})=0$, (see Remark~\ref{fixdef}).
\end{dfn}

 \begin{cor}[\cite{RRS,liu}]\label{lrrs}
 Let $T$ be a weakly mixing $\Z^2$-action
  on a Lebesgue probability space $(X,\mu)$, and let $L\subseteq\R^2$ be a $1$-dimensional direction. The following are equivalent
 \begin{enumerate}
 \item $T$ is direction $L$ weakly mixing. 
\item $T$ is direction $L$ ergodic,  \cite{RRS}. 
 \item $T$ has no direction $L$ Kronecker factors $K$,  {\rm (}i.e., $\ck_{\text{ap}_L}=\{\emptyset,X\}${\rm )}, \cite{liu}. 
 \item $\sigma^T$ has no direction $L^\perp$ walls: $\sigma^T(I)=0$ for any segment $I\subseteq\T^2$ with $I\perp L$,  \cite{RRS}.
\item \cite{liu}. Let $r>0$, $p=1,2$. For any $f,g\in L^2(X,\mu)$,
  \[\lim_{\ell\to\infty} 
  \frac{1}{\#(L^+(r,\ell))}
  \sum_{\vn\in L^+(r,\ell)}\lvert\langle U_T^\vn f,g\rangle-\langle f,1 \rangle\,\langle1,g\rangle\rvert^p=0,\]
 or equivalently, for any $A,B\in\cb$,
\[\lim_{\ell\to \infty} 
\frac{1}{\#(L^+(r,\ell))}
\sum_{\vn\in L^+(r,\ell)}|\mu(T^{-\vn}A\cap B)-\mu(A)\mu(B)|^p=0.\] 
 \item The Cartesian square $\Z^2$-action $T\times T$ is direction $L$ ergodic, or equivalently using (1) $\iff$ (2), $T\times T$ is direction $L$  
 weakly mixing, \cite{liu}.
 \item For any $r>0$, there exists a sequence $\{\vn_k\}_{k\in\N}\subseteq L^+(r)$
 so that 
  for any finite measurable partition $\xi$ on $(X,\mu)$, 
 the \emph{sequence entropy} satisfies 
  \[
 \limsup_{k\to\infty}\frac{1}{k}H_\mu\left(\,\bigvee_{j=1}^kT^{-\vn_j} \xi \right)=H_\mu(\xi),
 \]
 where $H_\mu(\xi)$ is the usual entropy of a partition {\rm (}see \cite{liu}{\rm )}.
 \item If  $T$ embeds into an $\R^2$-action $\overline T$, then the restriction, $\overline T^{|L}$, of $\overline T$ 
 to the direction $L$ (essentially an $\R$-action) is weakly mixing, \cite{RRS}. 

 \end{enumerate}
 \end{cor}
 
 \begin{cor}\label{rrsc}
 Let $T$ be an ergodic action of $\Z^2$ on a Lebesgue probability space $(X,\mu)$. 
 \begin{enumerate}
 \item If $T$ is not weakly mixing then it has no weakly mixing directions $L$, \cite{RRS}.
 \item If $T$ is weak mixing, the set $\cn$ of non weakly mixing directions  is {\rm (}at most{\rm )} countable. Moreover, for any countable set 
   $\cn$ of directions, there is a weakly mixing $\Z^2$-action $T$ with $\cn$ as its set of non weakly mixing directions, \cite{RRS}. 
 \item Direction $L$ weak mixing is a spectral property, and thus an isomorphism invariant, \cite{RRS}. 
 \item If $S$ is a factor of $T$, then $\cn_S\subseteq\cn_T$,  where $\cn_S$ and $\cn_T$ are the sets of non weakly mixing directions, 
 possibly with $\cn_T\backslash\cn_S\ne\emptyset$, \cite{RRS}.
 \item If $T$ is strongly mixing, then it is weakly mixing in every direction,  \cite{RRS}.
 \item The generic $\Z^2$-action $T$ is weakly mixing in every direction. Here, ``generic'' means residual in the weak 
 topology on the set of all $\Z^2$-actions $T$ on $(X,\mu)$, \cite{RRS}, {\rm (}see also \cite{kechris}{\rm )}.
   \item If $T$ is prime and not Kronecker then it is weakly mixing in every direction.
 \end{enumerate}
 \end{cor}
 
 Part (7) follows from 
 part (3) of Corollary~\ref{lrrs}, since in this case $T$ has no non-trivial factors at all. 
 
 \begin{rmk}
 In \cite{PR} it is shown that the \emph{Chacon $\Z^2$-action} $T$ 
 has \emph{minimal self joinings} (MSJ), and that this implies it is weakly mixing in every rational direction. 
Since 
 MSJ implies prime, part (7) of Corollary \ref{rrsc} 
implies that the Chacon $\Z^2$-action
$T$ is weakly mixing in all irrational directions as well, 
answering a question from \cite{RRS}.  
\end{rmk}

\section{The direction $L$  Kronecker  property}\label{old}

\begin{dfn}[\cite{liuxu}]
A $\Z^2$-action $K=\{K^\vn\}_{\vn\in\Z^2}$ on a Lebesgue probability space $(X,\mu)$, with Lebesgue $\sigma$-algebra $\cb$,
is called a \emph{direction $L$ Kronecker} $\Z^2$-action
if $\ck_{\text{ap}_L}=\cb$.
\end{dfn}

Using results in this paper, we obtain an alternative characterization of the direction $L$ Kronecker property.

\begin{prop}\label{kron}
A $\Z^2$-action $K$ is direction $L$ Kronecker if and only if 
its spectral measure satisfies $\sigma^K=\sigma^K_{\text{d}_L}$. 
\end{prop}

\noindent This follows from Theorem~\ref{main} and Proposition~\ref{JdLG}.

\smallskip

More generally, since $\ck_{\text{ap}_L}\subseteq\cb$ is a $T$-invariant sub $\sigma$-algebra,
and $(X,\mu)$ is a Lebesgue probability space,
$\ck_{\text{ap}_L}$ 
defines a factor $\Z^2$-action
$K$ of $T$, acting on a Lebesgue probability space $(Y,\nu)$.
In particular, if $\ca$ is the Lebesgue $\sigma$-algebra on $(Y,\nu)$, then $\ck_{\text{ap}_L}=\Psi^{-1}(\ca)$, 
where $\Psi:(X,\mu)\to (Y,\nu)$ is the factor map.
The factor $K$ has the 
direction $L$ Kronecker property, and is called 
the \emph{direction $L$ Kronecker factor} of $T$, \cite{liuxu}.
In particular, \cite{liuxu}, any $\Z^2$-action $T$ that is not direction $L$ weak mixing 
has a non-trivial direction $L$ Kronecker factor (see Definition~\ref{def1}).

Although many properties of direction $L$ Kronecker $\Z^2$-actions $K$ are described in \cite{liuxu}, 
no examples  are provided there. The goal of this section is to provide such examples.
The paper \cite{RRS} provides examples of $\Z^2$-actions $T$ with non weakly mixing directions $L$,
so by Proposition~\ref{kron}, these  have nontrivial direction $L$ Kronecker factors.     
Most of these examples in \cite{RRS} are produced using 
Gaussian actions $T$. But, before we discuss that method, we discuss some simpler examples.

\begin{example}[Kronecker actions]\label{alldir}
Suppose $K$ is an ergodic Kronecker  $\Z^2$-action 
on a Lebesgue probability space $(X,\mu)$, which means $\sigma^K=\sigma_{\text{d}}$.  
We claim that for every direction $L$, the action $K$ is direction $L$ Kronecker. 
To see this, fix a direction $L$.
Since $L^2_0(X,\mu)$ is separable, the set $\Sigma_K\subseteq\T^2$ of eigenvalues is countable, so that  
 $\sigma^K=\sum_{\vw\in \Sigma_K}a_\vw\delta_{\vw}$, $a_\vw>0$, is a countable sum. 
 Then there is a countable set $\Theta\subseteq[0,1)$ 
 so that $\Sigma_K\subseteq\bigcup_{\theta\in \Theta}L^\perp_{\theta\ve_1}$.
  It follows that $\sigma^K$ is concentrated on $\bigcup_{\theta\in \Theta}L^\perp_{\theta\ve_1}$, so that 
 $\sigma^K=\sigma_{\text{d}_L}$, and $K$ is direction $L$ Kronecker by Proposition~\ref{kron}.
 \end{example}

It is shown in \cite{liuxu} that if $K$ is a direction $L$ Kronecker $\Z^2$-action for at least two directions $L_1\ne L_2$, then it is Kronecker. 
It follows from Example~\ref{alldir} that such a $K$ is then directionally Kronecker in every direction.

 \begin{example}[Product type actions]\label{tensor}
 
 Suppose $T_1$ and $T_2$ are ergodic measure preserving transformations ($\Z$-actions)
 on $(X,\mu)$. Define the \emph{product type} $\Z^2$-action $T_1\otimes T_2$ 
 on $(X\times X,\cb\otimes\cb,\mu\times\mu)$ by 
 $(T_1\otimes T_2)^{(n_1,n_2)}(x,y)=(T_1^{n_1} x,T_2^{n_2} y)$, $\vn=(n_1,n_2)\in\Z^2$. 
 Note that $T_1\otimes T_2$ is ergodic, and is weakly mixing if and only if $T_1$ and $T_2$ are. 
 Clearly  $T_1\otimes I$ and $I\otimes T_2$ are both factors of $T_1\otimes T_2$.
 
 Let $\sigma^{T_1}$ and $\sigma^{T_2}$, on $\T$,  be the spectral measures
 for $T_1$ and $T_2$,  and let $\sigma^{T_1\otimes T_2}$, on $\T^2$,
be the spectral measure of the 
$\Z^2$-action $T_1\otimes T_2$, (as usual, in all of these we work orthogonal to the constant functions).
It is shown in \cite{RRS} that 
\[
\sigma^{T_1\otimes T_2}=
(\sigma^{T_1}\times \delta_0)+(\delta_0\times \sigma^{T_2})+(\sigma^{T_1}\times \sigma^{T_2}), 
\]
where 
$\sigma_{T_1}\times \delta_0$ is a wall 
measure concentrated on $L_1=\{(x,0):x\in\R\}\subseteq\T^2$ and 
$\delta_0\times \sigma_{T_2}$ is a wall measure concentrated on 
 $L_2=\{(0,y):y\in\R\}\subseteq\T^2$. 
 Thus vertical, $L_2=L_1^\perp$, and horizontal, $L_1=L_2^\perp$
are always non weakly mixing directions. 
 
 If $T_1$ and $T_2$ are both weakly mixing then
  $\sigma^{T_1}\times\sigma^{T_2}$ has no walls, 
  and thus vertical and horizontal directions are the only non-weakly mixing directions.
  It follows that $T_1\otimes T_2$  has exactly two directional Kronecker factors.
  It is easy to see that these are $T_1\otimes I$, which is the direction $L_2$ 
  Kronecker factor and $I\otimes T_2$, which is the direction $L_1$. 
In particular $T_1\otimes I$ and $I\otimes T_2$ are directionally Kronecker.
\end{example} 

Recall that an $\Z^2$ action $T$ is called (essentially) \emph{free} if $\mu\{x: \exists\ \vn\in\Z^2, \vn\ne{\vec 0}, \text{\rm\ so that } T^\vn x=x\}=0$. 
A standard result (see \cite{kt} or \cite{nadkarni}) is that for a free $\Z^2$-action,  the support of the spectral measure $\sigma^T$ is all of 
 $\T^2$ (i.e., $\sigma^T$ is concentrated on a dense subset of $\T^2$). 
 If $T$ is ergodic, the converse is also true (see the proof of Proposition~\ref{denses}). 
It is easy to see that neither $\Z^2$-action $T_1\otimes I$ or  $I\otimes T_2$ is free  (or ergodic).
For example, $(T_1\otimes I)^{(0,n)}(x,y)=(x,y)$ for all $n\in\Z$. And 
 indeed, the supports 
of their spectral measures ($L_2$ and $L_1$) are not all of  $\T^2$. 

\bigskip

Now we will construct  examples of Gaussian   $\Z^2$-actions (see Theorem~5.10 in \cite{RRS}) 
that are direction $L$ Kronecker for any given direction $L$. 
The benefit is that we can ensure these examples are free, but the price is the somewhat obscure nature of Gaussian examples. 

The general method of constructing a Gaussian $\Z^2$-action  starts with 
a finite positive Borel measure $\zeta$ on $\T^2$ that: (1) is 
\emph{non-atomic} and (2) is \emph{symmetric} in the sense that $\zeta(-E)=\zeta(E)$. 
The construction (called the ``Gaussian measure space (GMS) construction'' 
in \cite{RRS}) then produces a Lebesgue probability space
$(X,\mu)$ and weakly mixing $\Z^2$-action $T$ on $(X,\mu)$ satisfying 
\begin{equation}\label{exp}
\sigma^T=e^{\zeta}:=\zeta+(\zeta*\zeta)/{2!}+(\zeta*\zeta*\zeta)/{3!}+\dots\, .
\end{equation}
See \cite{CFS} or \cite{eli} for details in the $\Z$ case, the $\Z^2$ case being almost identical. 
We won't discuss these details here. 

\begin{example} If $L$ is a rational direction then $L$, projected 
to $\T^2$, is a circle through $0$, as is  
$L^\perp_{\vec 0}\subseteq\T^2$. Both are closed subgroups. 
Let $\zeta'$ be a nonatomic finite Borel measure on $\T^2$ supported 
on $L^\perp_{\vec 0}$ (e.g., such a $\zeta'$ can be obtained using the idea in Proposition~\ref{single}). 
Define a symmetric $\zeta$ by $\zeta(E)=(\zeta'(E)+\zeta'(-E))/2$. Let $K$ be the weakly mixing Gaussian 
$\Z^2$-action 
with $\sigma^K=e^\zeta$. In particular since $L^\perp_{\vec 0}$ is a subgroup, $\sigma^K$ is supported on $L^\perp_{\vec 0}$.
Thus by Proposition~\ref{kron}, $K$  
is a direction $L$ 
Kronecker action.
However,  like the previous examples,
$K$ is  not essentially  free, and the support of $\sigma^K$ is not dense.  
\end{example}

In order 
to produce free weakly mixing Gaussian directionally Kronecker 
$\Z^2$-actions in any direction $L$, we can use the following proposition. 

\begin{prop}\label{denses}
Suppose $\zeta$ is a 
nonatomic and symmetric finite positive  Borel measure on $\T^2$.
If $\text{supp}(\zeta)=\T^2$, then the 
Gaussian $\Z^2$-action $T$ with $\sigma^T=e^\zeta$ is weakly mixing and {\rm (}essentially{\rm )} free.
\end{prop}

\begin{proof}
Weak mixing implies ergodic, and $\zeta\ll\sigma^T$, so $\sigma^T$ has support $\T^2$. 

Now suppose $T$ is not essentially free. Then there exists $\vn_0\ne{\vec 0}$ such that $T^{\vn_0}={\rm id}$. 
The $\Z^2$-action $T$ can thus be viewed as an action $T_0$ of the quotient group 
$G=\Z^2/{\rm span}_\Z(\vn_0)\approx\Z\times \Z/g\Z$, where $g={\rm gcd}(a,b)$ for $\vn_0=(a,b)$.
Since $T$ is ergodic, $\sigma^T\sim\sigma^{T_0}$, because $T_0$ and $T$ are essentially the same action. On the 
other hand, ${\rm supp}(\sigma^{T_0})=\{\vs:\vn_0\cdot\vs\in\Z\}$ is a finite union of $g$ closed circles in $\T^2$ with slope $-a/b$.  
This contradicts the dense support of $\sigma^T$ in $\T^2$.
\end{proof}

\begin{example}
Let $L\subseteq\R^2$ be a direction in $\Z^2$. 
Enumerate a finite or 
countable set 
$\Theta:=\{\theta_1,\theta_2,\dots\}\subseteq\T$ 
so that 
$\cl:=\bigcup_k L^\perp_{\theta_k\ve_1}$ is dense in $\T^2$. 
One important case is when $L$ is a single irrational direction, $\Theta=\{0\}$, and $\cl=L^\perp_{\vec 0}$.
Note that the sets $L^\perp_{\theta_k\ve_1}$ are all cosets of the subgroup $L^\perp_{\vec 0}$ of $\T^2$. They are circles in the 
rational slope case.

For each $\theta_k\in \Theta$ let $\zeta_k$ be a nonatomic Borel probability measure 
on $L^\perp_{\theta_k\ve_1}$ whose support is  $L^\perp_{\theta_k\ve_1}$. Let 
$\zeta'=\sum_k\zeta_k/2^k$, and define $\zeta$ by $\zeta(E)=(\zeta'(E)+\zeta'(-E))/2$. 
Then $\zeta$ is a symmetric nonatomic Borel probability measure on $\T^2$ with support $\T^2$.
In particular, $\zeta$ is concentrated on a countable union of cosets $L^\perp_{\theta\ve_1}$. 
Let $K$ be the corresponding weakly mixing Gaussian $\Z^2$-action with $\sigma^K:=e^\zeta$.
By proposition~\ref{denses}, $K$ is essentially free. 

Note that if $\zeta_k$ and $\zeta_{k'}$ are concentrated, respectively, on the cosets $ L^\perp_{\theta_k\ve_1}$ and  $L^\perp_{\theta_{k'}\ve_1}$,
then the convolution
$\zeta_k*\zeta_{k'}$ is concentrated on their sum, the coset $L^\perp_{(\theta_k+\theta_{k'})\ve_1}$. 
 Letting $\Theta'=\text{span}_\Z\Theta$, which is countable, we have by \eqref{exp} that 
  $\text{supp}(e^\zeta)=\bigcup_{\theta\in \Theta'} L^\perp_{\theta\ve_1}$. Since $\sigma^K=e^\zeta$, Proposition~\ref{kron} implies 
  $K$ is a direction $L$ Kronecker $\Z^2$-action. 
\end{example}

\section{Concluding remarks}\label{conclusion}

Let $d>2$ and $1\le e<d$. An $e$-dimensional direction $L$ in $\Z^d$ is defined to be 
an $e$-dimensional subspace $L\subseteq\R^d$.  It is then easy to define ``spectral'' and ``dynamical''
direction $L$ weak mixing (as well as  direction $L$  ergodicity and Kronecker) for a $\Z^d$-action 
$T=\{T^\vn\}_{\vn\in\Z^d}$. But our proof that ``spectral'' and ``dynamical''
direction $L$ weak mixing are the same does not seem to easily generalize 
to this case.  In general, there are three possibilities for $L$: (i) $L$ is \emph{rational}, which means 
$L\cap \Z^d\approx \Z^e$, (ii)  $L$ is, in the terminology of \cite{RRS}, 
\emph{completely irrational}, which means $L\cap \Z^d=\{\vec 0\}$, and (iii) an  in between  case
$L\cap \Z^d\approx \Z^b$ for $0<b<e$. Only the first two cases occur when $d=2$ and $e=1$, and this fact is used for 
our proof of Theorem~\ref{main}. But we conjecture that Theorem~\ref{main} is true in general.

\section*{Acknowledgement} While completing work on this paper the authors received the preprint \cite{liuetal} which, as well as proving other interesting results, 
proves that ``dynamical'' and ``spectral'' direction $L$ weak mixing are the same for $\Z^2$-actions $T$. The proof in  \cite{liuetal}
is based on a completely different idea than the proof in this paper. 

\bigskip

\section*{Declaration of Generative AI in the writing process} The overall structure of this paper and all major arguments are due to the authors.
For some proofs, AI (Anthropic Claude Sonnet 5.5) was used to suggest paths forward, sometimes having first been given a 
proof strategy or an outline. In many cases, 
producing a correct argument in this way required a substantial back and forth collaboration between the authors and AI.
AI was also employed for proofreading, but in the final document, no AI generated text was used, and   
the authors take full responsibility for all the content.

\end{document}